\documentclass[11pt]{article}
\usepackage{epsfig}

\def\be{\begin{equation}}
\def\ee{\end{equation}}
\def\bea{\begin{eqnarray}}
\def\eea{\end{eqnarray}}
\def\bes{\begin{eqnarray*}}
\def\ees{\end{eqnarray*}}

\def\nn{\nonumber}
\def\lb{\label}
\def\bs{\setminus}

\def\R{{\bf R}}
\def\C{{\bf C}}
\def\Z{{\bf Z}}

\def\N{{\bf N}}
\def\U{{\bf U}}

\def\Q{{\bf Q}}

\def\CP{{\bf CP}}
\def\HP{{\bf HP}}
\def\CaP{{\bf CaP}}
\def\RP{{\bf RP}}

\def\aa{{\alpha}}
\def\bb{{\beta}}
\def\ga{{\gamma}}

\def\th{{\theta}}

\def\om{{\omega}}
\def\Om{{\Omega}}
\def\ep{{\epsilon}}
\def\lm{{\lambda}}
\def\Lm{{\Lambda}}

\def\sg{{\sigma}}
\def\dm{{\diamond}}

\def\<{{\langle}}
\def\>{{\rangle}}

\def\P{{\cal P}}

\def\diag{{\rm diag}}

\def\rank{{\rm rank}}

\def\Sp{{\rm Sp}}

\def\mod{{\rm mod}}
\def\ol{\overline}

\def\hb{\vrule height0.18cm width0.14cm $\,$}

\title{Elliptic closed geodesics on non-negatively curved  Finsler manifolds}

\author{Wei Wang\thanks{Partially supported by NSFC No. 12025101.
E-mail: wangwei@math.pku.edu.cn  }\\Key Laboratory of Pure and Applied Mathematics\\
School of Mathematical Sciences \\ Peking University, Beijing 100871 \\
PEOPLES REPUBLIC OF CHINA \\ }
\begin{document}
\maketitle

\begin{abstract}
{\it In this paper, we prove that if a compact simply-connected
manifold $(M,F)$ with a bumpy irreversible Finsler metric $F$ whose
flag curvature $K$ satisfies $K\ge 0$ and $H^*(M;\Q)\cong
T_{d,n+1}(x)$ for some $d\in2\N$ and $n\in\N$ carries exactly
$\frac{dn(n+1)}{2}$ (the number of closed geodesics in Katok's
example) closed geodesics, then there are at least $2$ irrationally
elliptic closed geodesics on $(M, F)$.  }
\end{abstract}



\renewcommand{\theequation}{\thesection.\arabic{equation}}
\renewcommand{\thefigure}{\thesection.\arabic{figure}}

\setcounter{figure}{0}
\setcounter{equation}{0}
\section{Introduction and main results}

There is a famous conjecture in Riemannian geometry which claims the
existence of infinitely many distinct closed geodesics on every
compact Riemannian manifold $M$. This conjecture has been proved for
many cases (cf. for example \cite{BTZ1}, \cite{BTZ2}, \cite{Ban1}),
specially Gromoll and Meyer proved the conjecture in \cite{GrM} of
1969 for any compact $M$ provided the Betti number sequence
$\{b_p(\Lm M)\}_{p\in\N}$ of the free loop space $\Lm M$ of $M$ is
unbounded. Then in \cite{ViS1} of 1976, for compact simply-connected
manifold $M$, Vigu\'e-Poirrier and Sullivan proved this Betti number
sequence is bounded if and only if $M$ satisfies \be H^*(M;\Q)\cong
T_{d,n+1}(x)=\Q[x]/(x^{n+1}=0)   \lb{1.1}\ee with a generator $x$ of
degree $d\ge 2$ and height $n+1\ge 2$, where $\dim M=dn$. Note that
when $d$ is odd, then $x^2=0$ and $n=1$, or when $n=1$, $M$ is
rationally homotopic to $S^d$ (cf. Remark 2.5 of \cite{Rad1} and
\cite{Hin}). Among these manifolds, only for Riemannian $S^2$ this
conjecture was proved by the works \cite{Fra} of Franks in 1992 and
\cite{Ban2} of Bangert in 1993. When Finsler manifolds are
considered, the situation changes dramatically. It was quite
surprising that Katok found some irreversible Finsler metrics on compact
rank one symmetric spaces which possess finitely many distinct
closed geodesics in \cite{Kat} of 1973. There Katok carried out also
detailed constructions on the Finsler metric on $S^d$ which possess
precisely $2[(d+1)/2]$ distinct closed geodesics and all of which
are non-degenerate and elliptic. The geometry of Katok's metrics was
further studied by Ziller in \cite{Zil} of 1982. He constructed in
detail Finsler metrics for complex projective spaces $\CP^n$ (with
$d=2$), quaternionic projective spaces $\HP^n$ (with $d=4$), and the
Cayley plane $\CaP^2$ (with $d=8$ and $n=2$) which possess precisely
$n(n+1)$, $2n(n+1)$, and $24$ distinct closed geodesics
respectively, i.e., precisely $\frac{dn(n+1)}{2}$ distinct closed
geodesics. Based on Katok's work, Anosov in \cite{Ano} conjectured in ICM of 1974
that the lower bound of the number of distinct closed geodesics on
any Finsler sphere $(S^d, F)$ should be $2[\frac{d+1}{2}]$, i.e.,
the number of closed geodesics in Katok's example. Thus it is
natural to generalize this conjecture to all compact simply-connected Finsler manifolds satisfying (\ref{1.1}), i.e., the lower
bound of the number of distinct closed geodesics on such manifolds
is $\frac{dn(n+1)}{2}$.

The following result of Duan,  Long and Wang shows that this
conjecture is true under some additional conditions:

{\bf Theorem 1.1.} (cf. Theorem 1.2 of \cite{DLW2}) {\it For any
compact simply-connected manifold $(M,F)$ with a bumpy irreversible
Finsler metric $F$ and $H^*(M;\Q)\cong T_{d,n+1}(x)$ for some even
integer $d\ge 2$ and some integer $n\ge 1$, if the flag curvature
$K\ge 0$, then there exist at least $\frac{dn(n+1)}{2}$ distinct
prime closed geodesics.}

Recall that a closed curve on a Finsler manifold is a closed geodesic if it is locally the shortest
path connecting any two nearby points on this curve (cf. \cite{BCS}, \cite{She}). As usual, on any Finsler
manifold $(M,F)$, a closed geodesic $c: S^1=\R/\Z\to M$ is {\it prime} if it is not a multiple
covering (i.e., iteration) of any other closed geodesics. Here the $m$-th iteration $c^m$ of $c$ is
defined by $c^m(t)=c(mt)$. The inverse curve $c^{-1}$ of $c$ is defined by $c^{-1}(t)=c(1-t)$ for
$t\in \R$. Note that unlike the Riemannian case, the inverse curve $c^{-1}$ of a closed geodesic $c$
on an irreversible Finsler manifold need not be a geodesic. We call two prime closed geodesics $c$
and $d$ {\it distinct} on a Finsler manifold, if there exists no $\th\in (0,1)$ such that
$c(t)=d(t+\th)$ for all $t\in\R$. On a reversible Finsler
(or Riemannian) manifold, i.e., $F(X)=F(-X)$ for
$\forall X\in TM$, two closed
geodesics $c$ and $d$ are called { \it geometrically distinct} if $c(S^1)\neq d(S^1)$, i.e., their
image sets in $M$ are distinct. We shall omit the word {\it distinct} when we talk about more than
one prime closed geodesic for simplicity.

For a closed geodesic $c$ on the $n$-dimensional manifold $(M,\,F)$, denote by $P_c$ the linearized
Poincar\'{e} map of $c$. Then $P_c\in \Sp(2n-2)$ is a symplectic matrix. As usual, for any $M\in \Sp(2k)$, we
define the {\it elliptic height} $e(M)$ of $M$ to be the total algebraic multiplicity of all
eigenvalues of $M$ on the unit circle $\U=\{z\in\C|\; |z|=1\}$ in the complex plane $\C$. Since $M$
is symplectic, $e(M)$ is even and $0\le e(M)\le 2k$. A closed geodesic $c$ is {\it elliptic} if
$e(P_c)=2(n-1)$, i.e., all the eigenvalues of $P_c$ locate on $\U$; {\it hyperbolic} if $e(P_c)=0$,
i.e., all the eigenvalues of $P_c$ locate away from $\U$; {\it non-degenerate} if $1$ is not an
eigenvalue of $P_c$; {\it
irrationally elliptic} if there exists $P\in \Sp(2n-2)$ such that
$P^{-1}P_cP=\diag(R(\theta_1), R(\theta_2),\ldots, R(\theta_{(n-1)}))$,
where $R(\theta_i)=\left(\matrix{\cos\th_i & -\sin\th_i\cr
                      \sin\th_i &  \cos\th_i\cr}\right)$
with rotation angles $\th_i$ satisfying $\frac{\theta_i}{\pi}\notin\Q$ for
$1\le i\le n-1$. A Finsler manifold $(M,\,F)$ is called {\it bumpy} if all the closed geodesics
and their iterates on $M$ are non-degenerate (cf. \cite{Abr}).

Due to  Katok's example, one may tend to believe that all the
geodesics on $(M, F)$ with $H^*(M;\Q)\cong T_{d,n+1}(x)$ should be
non-degenerate and elliptic provided the number of closed geodesics
on $(M, F)$ is finite. The following is an attempt to this problem:

{\bf Theorem 1.2.} {\it If a compact simply-connected manifold
$(M,F)$ with a bumpy irreversible Finsler metric $F$ whose flag curvature $K$
satisfies $K\ge 0$ and $H^*(M;\Q)\cong T_{d,n+1}(x)$ for some
$d\in2\N$ and $n\in\N$ carries exactly $\frac{dn(n+1)}{2}$ (the number
of closed geodesics in Katok's example) closed geodesics, then
there are at least $2$ irrationally elliptic closed geodesics on
$(M, F)$. }

{\bf Remark 1.3.} Note that for the standard Riemannian metric on
$S^n$ with constant curvature $1$ and the Fubini-Study metric on
$\CP^n$, there are no irrationally elliptic closed geodesics, cf. \cite{Bes}.

From 2005, the Maslov-type index theory for symplectic paths has been applied to study the closed
geodesic problem. In 2005, Bangert and Long proved the existence of at least two distinct closed
geodesics on every Finsler $(S^2,F)$ (which was published as \cite{BaL} in 2010). Since then in the
last twenty years, a great number of results on the multiplicity and stability of closed geodesics on
Finsler manifolds has appeared, for which we refer readers to
\cite{D1, D2}, \cite{DL1, DL2}, \cite{DuL1, DuL2, DuL3},
\cite{DQ}, \cite{DX}, \cite{GG}, \cite{GGM},  \cite{GM},
\cite{HiR}, \cite{Lon4}, \cite{LoD},
\cite{LX}, \cite{LLX}, \cite{LoW}, \cite{Rad4, Rad5, Rad6},
\cite{Wan1, Wan2, Wan3}, \cite{XiL}
and the references therein.

In 2008, Long and Wang in \cite{LoW} proved there exist
at least two irrationally elliptic closed geodesics
on $(S^2, F)$ provided the number of closed geodesics on $(S^2, F)$
is finite. In 2012, Wang in \cite{Wan2} proved the existence of at least $2[\frac{n+1}{2}]$ distinct closed
geodesics on a bumpy Finsler $(S^n,F)$ when the flag curvature $K$ satisfies
$\frac{\lm^2}{(1+\lm)^2}<K\le 1$, where the reversibility $\lambda=\lambda(M,F)$ was introduced
by Rademacher in \cite{Rad3} as
\be   \lambda=\max\{F(-X)\ |\ X\in TM,\ F(X)=1\}\ge 1.   \lb{1.2}\ee
Note that the lower bound $2[\frac{n+1}{2}]$ on the number of distinct closed geodesics for spheres $S^n$
is sharp due to the above mentioned examples of Katok. Besides spheres, in \cite{Rad4}, Rademacher obtained some
multiplicity and stability results of closed geodesics on compact simply-connected manifolds satisfying
some pinching conditions. For example, it was proved
that a Finsler $(\CP^n,F)$ with $n\ge 7$ and a bumpy Finsler metric possessing only finitely many distinct
closed geodesics and satisfying $\left(\frac{2}{n+1}\frac{\lambda}{1+\lambda}\right)^2\le K\le 1$, carries
always at least $2n$ distinct closed geodesics, and at least $(n-3)$ of them are non-hyperbolic, where
$\lambda$ and $K$ denote the reversibility and the flag curvature of $(\CP^n, F)$ respectively. In
\cite{Rad5}, Rademacher further proved the existence of two prime closed geodesics on any $\CP^2$ with a bumpy irreversible Finsler metric. In \cite{DLW1},  Duan,
 Long and  Wang proved the existence of at
least two prime closed geodesics on $(M,F)$ with $H^*(M;\Q)\cong T_{d,n+1}(x)$ and a bumpy irreversible
Finsler metric $F$.

In this paper, we denote by $\N$, $\N_0$, $\Z$, $\Q$, $\R$, and $\C$ the sets of natural integers,
non-negative integers, integers, rational numbers, real numbers, and complex numbers respectively.
We use only singular homology modules with $\Q$-coefficients. We use the functions
\be \left\{\matrix{[a]=\max\{k\in\Z\,|\,k\le a\}, &
           E(a)=\min\{k\in\Z\,|\,k\ge a\} , \cr
    \varphi(a)=E(a)-[a], &\{a\}=a-[a].  \cr}\right. \lb{1.3}\ee
Especially, $\varphi(a)=0$ if $ a\in\Z\,$, and $\varphi(a)=1$ if $a\notin\Z\,$.

\setcounter{figure}{0}
\setcounter{equation}{0}
\section{Morse theory of closed geodesics}

Let $M=(M,F)$ be a compact Finsler manifold, the space $\Lambda=\Lambda M$ of $H^1$-maps
$\gamma:S^1\rightarrow M$ has a natural structure of Riemannian Hilbert manifolds on which the
group $S^1=\R/\Z$ acts continuously by isometries (cf. \cite{Kli}). This action is defined by
$(s\cdot\ga)(t)=\ga(t+s)$ for all $\gamma\in\Lm$ and $s, t\in S^1$. For any $\gamma\in\Lambda$,
the energy functional is defined by
\be E(\gamma)=\frac{1}{2}\int_{S^1}F(\gamma(t),\dot{\gamma}(t))^2dt.  \lb{2.1}\ee
It is $C^{1,1}$ and invariant under the $S^1$-action, cf. \cite{Mer}. The critical points of $E$ of positive
energies are precisely the closed geodesics $c: S^1\to M$. The index form of the functional
$E$ is well defined on any closed geodesic $c$ on $M$, which we denote by $E''(c)$. As usual,
we denote by $i(c)$ and $\nu(c)$ the Morse index and nullity of $E$
at $c$ (cf. \cite{Kli}, \cite{Mor} and \cite{Rad1}). In the following, we denote by
\be \Lm^\kappa=\{d\in \Lm\;|\;E(d)\le\kappa\},\quad \Lm^{\kappa-}=\{d\in \Lm\;|\; E(d)<\kappa\},
  \quad \forall \kappa\ge 0. \lb{2.2}\ee
For a closed geodesic $c$ we set $ \Lm(c)=\{\ga\in\Lm\mid E(\ga)<E(c)\}$.

For $m\in\N$ we denote the $m$-fold iteration map
$\phi_m:\Lambda\rightarrow\Lambda$ by $\phi_m(\ga)(t)=\ga(mt)$, for all
$\,\ga\in\Lm, t\in S^1$, as well as $\ga^m=\phi_m(\gamma)$. If $\gamma\in\Lambda$
is not constant then the multiplicity $m(\gamma)$ of $\gamma$ is the order of the
isotropy group $\{s\in S^1\mid s\cdot\gamma=\gamma\}$. For a closed geodesic $c$,
the mean index $\hat{i}(c)$ is defined as usual by
$\hat{i}(c)=\lim_{m\to\infty}i(c^m)/m$. Using singular homology with rational
coefficients we consider the following critical $\Q$-module of a closed geodesic
$c\in\Lambda$:
\be \overline{C}_*(E,c)
   = H_*\left((\Lm(c)\cup S^1\cdot c)/S^1,\Lm(c)/S^1\right). \lb{2.3}\ee

{\bf Lemma 2.1.} (cf. Satz 6.11 of \cite{Rad2} ) {\it Let $c$ be a
prime closed geodesic on a bumpy Finsler manifold $(M,F)$. Then there holds}
$$ \overline{C}_q( E,c^m) = \left\{\matrix{
     \Q, &\quad {\it if}\;\; i(c^m)-i(c)\in 2\Z\;\;{\it and}\;\;
                   q=i(c^m),\;  \cr
     0, &\quad {\it otherwise}. \cr}\right.  $$

{\bf Definition 2.2.} (cf. Definition 1.6 of \cite{Rad1}) {\it For a
closed geodesic $c$, let $\ga_c\in\{\pm\frac{1}{2},\pm1\}$ be the
invariant defined by $\ga_c>0$ if and only if $i(c)$ is even, and
$|\ga_c|=1$ if and only if $i(c^2)-i(c)$ is even. }

{\bf Theorem 2.3.} (cf. Theorem 3.1 of \cite{Rad1} and Satz 7.9 of
\cite{Rad2}) {\it Let $(M,F)$ be a compact simply-connected bumpy
Finsler manifold with $\,H^{\ast}(M,\Q)=T_{d,n+1}(x)$ for $d\in 2\N$
and $n\in\N$. Denote the  prime closed geodesics on $(M,F)$ with positive
mean indices by $\{c_j\}_{1\le j\le q}$ for some $q\in\N$. Then the
following identity holds \be
\sum_{j=1}^q\frac{\ga_{c_j}}{\hat{i}(c_j)}=B(d,n) =
     -\frac{n(n+1)d}{2d(n+1)-4}.   \lb{2.4}\ee
where $\dim M=dn$.}

\medskip

{\bf Lemma 2.4.} (cf. Theorem 2.4 of \cite{Rad1}) {\it Let $M$ be a
compact simply-connected
manifold with $H^*(M;\Q)\cong T_{d,n+1}(x)$ for some integers $d\in
2\N$ and $n\in\N$.
Let $D=d(n+1)-2$ and
\bea \Om(d,n) = \{k\in 2\N-1&\,|\,& iD\le k-(d-1)=iD+jd\le iD+(n-1)d\;  \nn\\
         && \mbox{for some}\;i\in\N\;\mbox{and}\;j\in \{0,1,\ldots,n-1\}\}. \nn\eea
Then the Betti numbers of the free loop space of $M$ are given by
\bea  b_p(\Lm M)
&=& \rank H_p(\Lm M/S^1,\Lm^0 M/S^1;\Q)\nn\\
&=& \left\{\matrix{
    0, & \quad \mbox{if}\ p\ \mbox{is even or}\ p\le d-2,  \cr
    \left[\frac{p-(d-1)}{d}\right]+1, & \quad \mbox{if}\ p\in 2\N-1\;\mbox{and}\;d-1\le p < d-1+(n-1)d, \cr
    n+1, & \quad \mbox{if}\ p\in \Om(d,n), \cr
    n, & \quad \mbox{otherwise}. \cr}\right.\lb{2.5}\eea
}

{\bf Proof.} For a compact and simply-connected Finsler manifold
$(M, F)$ with $H^\ast(M;\Q)\cong T_{d, n+1}(x)$ for $d\in2\N$ and
$n\in\N$, the Poincar\'e series of its free loop space was given by
Theorem 2.4 of \cite{Rad1} \bea \sum_{q=0}^{+\infty}b_qt^q&=&P(\Lm
M/S^1,\Lm^0 M/S^1)(t)
=t^{d-1}\left(\frac{1}{1-t^2}+\frac{t^{d(n+1)-2}}{1-t^{d(n+1)-2}}\right)\frac{1-t^{dn}}{1-t^d}\nn\\
&=&t^{d-1}\left(\frac{1}{1-t^2}+\frac{t^D}{1-t^D}\right)\frac{1-t^{dn}}{1-t^d}\nn\\
&=&t^{d-1}\left(\sum_{i=0}^{+\infty}t^{2i}
+\sum_{i=1}^{+\infty}t^{iD}\right)\left(\sum_{j=0}^{n-1}t^{jd}\right).
\lb{2.6}\eea
Since $d\in2\N$, we have $b_q=0$ for $q\in2\N$. Clearly,
$b_q=0$ for $q\le d-2$.

If $p\in 2\N-1$ and $d-1\le p < d-1+(n-1)d=nd-1$.
Note that $D=d(n+1)-2=nd+d-2\ge nd$ since $d\ge 2$,
therefore $b_p$ counts the number of $t^p$'s
in the expression
\bea t^{d-1}\left(\sum_{i=0}^{+\infty}t^{2i}
\right)\left(\sum_{j=0}^{n-1}t^{jd}\right)\lb{2.7}
\eea
Write $p=d-1+kd+2l$ for some integers $k, l$ satisfying
$0\le k \le n-2$
and $0\le 2l\le d-2$, then the  $t^p$'s are given by
$(i, j)\in\{(\frac{kd}{2}+l, 0), (\frac{(k-1)d}{2}+l, 1),\ldots,
(l, k)\}$. Hence $b_q=k+1=[\frac{p-(d-1)}{d}]+1$ as desired.

Now consider $p\in 2\N-1$ and $p\ge nd-1$.  Write $p=d-1+kd+2l$ for
some integers $k, l$ satisfying $k \ge n-1$ and $0\le 2l\le d-2$,
then the  $t^p$'s in (\ref{2.7})  are given by $(i,
j)\in\{(\frac{kd}{2}+l, 0), (\frac{(k-1)d}{2}+l, 1),\ldots,
(\frac{(k-(n-1))d}{2}+l, n-1)\}$. Hence the number of $t^p$'s in
(\ref{2.7}) is $n$.

Now we consider the number of $t^p$'s in the expression
 \bea t^{d-1}\left(
\sum_{i=1}^{+\infty}t^{iD}\right)\left(\sum_{j=0}^{n-1}t^{jd}\right)
\lb{2.8}\eea If $p\in\Omega(d, n)$, write $p=(d-1)+i_\ast D+j_\ast d$
 for some $i_\ast\in\N$
and $j_\ast\in \{0,1,\ldots,n-1\}$, then there is exactly one term
$(i, j)=(i_\ast, j_\ast)$ contributes to $t^p$'s  in (\ref{2.8}).
In fact, this follows from $iD>(i-1)D+(n-1)d$ for
$i\in\N$ since $d\ge 2$. Hence $b_p=n+1$ for $p\in\Omega(d, n)$ and $b_p=n$ otherwise.
 \hfill\hb

\medskip

{\bf Theorem 2.5.} (cf. Theorem I.4.3 of \cite{Cha}) {\it Suppose that there exist only finitely many
prime closed geodesics $\{c_k\}_{1\le k\le q}$ on a Finsler manifold $(M,F)$. Define
\be M_p=\sum_{1\le k\le q,\ m\ge 1}\dim \ol{C}_p(E,c_k^m),\qquad\forall\ p\in\Z. \lb{2.9}\ee
Then the following Morse inequalities holds for every $p\in \N_0$,
\bea M_p &\ge& b_p(\Lm M),\lb{2.10}\\
M_p - M_{p-1} + \cdots +(-1)^{p}M_0
&\ge& b_p(\Lm M) - b_{p-1}(\Lm M)+ \cdots + (-1)^{p}b_0(\Lm M).  \lb{2.11}\eea}

\setcounter{figure}{0} \setcounter{equation}{0}
\section{Index iteration theory for symplectic paths}

In this section, we  review the index iteration theory developed by Y.  Long
and his coworkers, which will be used in the next section. In \cite{Lon1} of 1999, Y. Long established the
basic normal form decomposition of symplectic matrices. Based on
this result he further established the precise iteration formulae of
indices of symplectic paths in \cite{Lon2} of 2000.

As in \cite{Lon2}, denote by \bea N_1(\lm, b) &=&
\left(\matrix{\lm
& b\cr
0 & \lm\cr}\right), \qquad {\rm for\;}\lm=\pm 1, \; b=\pm 1, 0, \lb{3.1}\\
R(\th) &=& \left(\matrix{\cos\th & -\sin\th \cr
                           \sin\th & \cos\th\cr}\right), \qquad {\rm for\;}\th\in (0,\pi)\cup (\pi,2\pi), \lb{3.2}\\
N_2(e^{\sqrt{-1}\th}, b) &=& \left(\matrix{ R(\th) & b \cr
                  0 & R(\th)\cr}\right), \qquad {\rm for\;}\th\in (0,\pi)\cup (\pi,2\pi)\;\; {\rm and}\; \nn\\
        && \qquad b=\left(\matrix{b_1 & b_2\cr
                                  b_3 & b_4\cr}\right)\; {\rm with}\; b_j\in\R, \;\;
                                         {\rm and}\;\;
                                         b_2\not= b_3, \lb{3.3}\\
D(\lm) &=& \left(\matrix{\lm & 0\cr
0 & \lm^{-1}\cr}\right), \qquad {\rm for\;}\lm=\pm 2. \lb{3.4}
\eea
Here $N_2(e^{\sqrt{-1}\th}, b)$ is non-trivial if
$(b_2-b_3)\sin\theta<0$, and trivial if
$(b_2-b_3)\sin\theta>0$.

As in \cite{Lon2}, the $\diamond$-product of any two block-wise real
matrices is defined by
$$ \left(\matrix{A_1 & B_1\cr
                 C_1 & D_1\cr}\right)_{2i\times 2i}\diamond
   \left(\matrix{A_2 & B_2\cr
                 C_2 & D_2\cr}\right)_{2j\times 2j}
=\left(\matrix{A_1 & 0 & B_1 & 0  \cr
                 0 & A_2 & 0 & B_2\cr
               C_1 & 0 & D_1 & 0  \cr
                 0 & C_2 & 0 & D_2}\right). $$

For any $M\in\Sp(2n)$, the homotopy set $\Omega(M)$ of $M$ in
$\Sp(2n)$ is defined by
$$ \Om(M)=\{N\in\Sp(2n)\,|\,\sg(N)\cap\U=\sg(M)\cap\U\;\mbox{and}
                    \;\nu_{\om}(N)=\nu_{\om}(M),\, \forall\om\in\U\}, $$
where we denote by $\sg(M)$ the spectrum of $M$, by $\Om^0(M)$ the
path connected component of $\Om(M)$ containing $M$, and
$\nu_{\om}(M)\equiv\dim_{\C}\ker_{\C}(M-\om I)$ for $\om\in\U$.

{\bf Lemma 3.1.} (cf. \cite{Lon1}, Lemmas 9.1.5, 9.1.6 and List 9.1.12 of
\cite{Lon3}) {\it For $M\in\Sp(2n)$ and $\om\in\U$, the splitting
number $S_M^\pm(\om)$ (cf. Definition 9.1.4 of \cite{Lon3})
satisfies $S_N^{\pm}(\om)=S_M^{\pm}(\om)$ for $N\in\Omega^0(M)$ and
\bea
S_M^{\pm}(\om) &=& 0, \qquad {\it if}\;\;\om\not\in\sg(M).  \lb{3.5}\\
S_M^+(\om) &=& S_M^-(\overline \om). \lb{3.6}\eea

For any $M_i\in\Sp(2n_i)$ with $i=0$ and $1$, there holds \be
S^{\pm}_{M_0\diamond M_1}(\om) = S^{\pm}_{M_0}(\om) +
S^{\pm}_{M_1}(\om),
    \qquad \forall\;\om\in\U. \lb{3.7}\ee

(i) $(S_M^+(1),\,S_M^-(1))=(1, \,1)$ for $M=N_1(1,\,b)$ with $b=0$
or $1$.

(ii) $(S_M^+(1),\,S_M^-(1))=(0, \,0)$ for $M=N_1(1,\,-1)$.

(iii) $(S_M^+(-1),\,S_M^-(-1))=(1, \,1)$ for $M=N_1(-1,\,b)$ with
$b=0$ or $-1$.

(iv) $(S_M^+(-1),\,S_M^-(-1))=(0, \,0)$ for $M=N_1(-1,\,1)$.

(v) $(S_M^+(e^{\sqrt{-1}\theta}),\,S_M^-(e^{\sqrt{-1}\theta}))=(0,
\,1)$ for $M=R(\theta)$ with $\theta\in(0,\,\pi)\cup(\pi,\,2\pi)$.

(vi) $(S_M^+(\om),\,S_M^-(\om))=(1, \,1)$ for $M=N_2(\om,\,b)$ being
non-trivial with $\om=e^{\sqrt{-1}\theta}\in\U\setminus\R$.

(vii) $(S_M^+(\om),\,S_M^-(\om))=(0, \,0)$ for $M=N_2(\om,\,b)$
being trivial with $\om=e^{\sqrt{-1}\theta}\in\U\setminus\R$.

(viii) $(S_M^+(\om),\,S_M^-(\om))=(0, \,0)$ for any $\om\in\U$ and
$M\in\Sp(2n)$ satisfying $\sigma(M)\cap\U=\emptyset$.}

We have the following normal form decomposition theorem for symplectic
matrix.

{\bf Theorem 3.2.} (cf. \cite{Lon2} and Theorem 1.8.10 of
\cite{Lon3}) {\it For any $M\in\Sp(2n)$, there exists a path $f\in
C([0,1],\Om^0(M))$ such that $f(0)=M$ and \be f(1) =
M_1\diamond\cdots\diamond M_k,  \lb{3.8}\ee where each
$M_i$ is a basic normal form listed in (\ref{3.1})-(\ref{3.4}) for
$1\leq i\leq k$.}

\medskip
For every $\ga\in\mathcal{P}_\tau(2n)\equiv\{\ga\in
C([0,\tau],Sp(2n))\ |\ \ga(0)=I_{2n}\}$, we extend $\ga(t)$ to $t\in
[0,m\tau]$ for every $m\in\N$ by \be
\ga^m(t)=\ga(t-j\tau)\ga(\tau)^j,\qquad \forall\;j\tau\le t\le
(j+1)\tau \;\;
               {\rm and}\;\;j=0, 1, \ldots, m-1, \lb{3.9}\ee
as in P.114 of \cite{Lon1}. As in \cite{LoZ} and \cite{Lon3}, we
denote the Maslov-type indices of $\ga^m$ by
$(i(\ga,m),\nu(\ga,m))$. Note that by \cite{Liu} and \cite{LLo},
the Morse indices of an oriented closed geodesic $c$ coincide to the
Maslov-type indices of the fundamental solution $\ga_c$ of the
linearized Hamiltonian system at $c$, i.e., $(i(c^m),\nu(c^m)) =
(i(\ga_c,m),\nu(\ga_c,m))$.

Denote by \bea
\Sp_{cnu}(2n) &=& \{M\in \Sp(2n)\;|\; \dim(M^m-I)=\dim(M-I),\;\forall\;m\in\N\},  \nn\\
\Sp_{vnu}(2n) &=& \Sp(2n)\bs \Sp_{cnu}(2n).  \nn\eea Note that every
$M\in \Sp_{vnu}(2n)$ must possess at least one eigenvalue
$e^{\sqrt{-1}\th}\in\sg(M)$ with $\th\in (0,2\pi)\cap\pi\Q$. Then
for every $M\in\Sp(2n)$ we define \be  \check{m}(M) =
\left\{\matrix{
   & \min\{k\in\N\;|\;k\th\in 2\pi\N,\;e^{\sqrt{-1}\th}\in\sg(M)\;{\rm with}\;\th\in (0,2\pi)\cap \pi\Q\},
            &{\rm if}\;M\in\Sp_{vnu}(2n), \cr
   & +\infty, & {\rm if}\;M\in\Sp_{cnu}(2n). \cr}\right.  \lb{3.10}\ee
Note that $\check{m}(M)\ge 2$ holds for every $M\in\Sp(2n)$. For a
closed geodesic $c$ on a Finsler manifold $(M,F)$, using the
linearized Poincar\'e map $P_c$ we define $\check{m}(c) =
\check{m}(P_c)$ correspondingly. When $(M,F)$ is bumpy, then
$\check{m}(c)=+\infty$ holds for every closed geodesic $c$ on $(M,F)$.

The following is the precise index iteration formulae for symplectic
paths, which is due to Y. Long (cf. Theorem 8.3.1 and Corollary
8.3.2 of \cite{Lon3}).

{\bf Theorem 3.3.} {\it Let $\gamma\in\P_\tau(2n)$, then there
exists a path $f\in C([0,1],\Omega^0(\gamma(\tau)))$ such that
$f(0)=\gamma(\tau)$ and \bea f(1)=&&N_1(1,1)^{\diamond p_-} \diamond
I_{2p_0}\diamond N_1(1,-1)^{\diamond p_+} \diamond
N_1(-1,1)^{\diamond q_-} \diamond (-I_{2q_0})\diamond
N_1(-1,-1)^{\diamond q_+}\nn\\
&&\diamond R(\theta_1)\diamond\cdots\diamond R(\theta_r)
\diamond N_2(\omega_1, u_1)\diamond\cdots\diamond N_2(\omega_{r_*}, u_{r_*}) \nn\\
&&\diamond N_2(\lm_1, v_1)\diamond\cdots\diamond N_2(\lm_{r_0},
v_{r_0}) \diamond M_0 \lb{3.11}\eea where $N_2(\omega_j, u_j)$'s
are non-trivial and   $N_2(\lm_j, v_j)$'s  are trivial basic normal
forms; $\sigma (M_0)\cap U=\emptyset$; $p_-$, $p_0$, $p_+$, $q_-$,
$q_0$, $q_+$, $r$, $r_*$ and $r_0$ are non-negative integers;
$\omega_j=e^{\sqrt{-1}\alpha_j}$, $ \lambda_j=e^{\sqrt{-1}\beta_j}$;
$\theta_j$, $\alpha_j$, $\beta_j$ $\in (0, \pi)\cup (\pi, 2\pi)$;
these integers and real numbers are uniquely determined by
$\gamma(\tau)$. Then using the functions defined in (\ref{1.3}) \bea
i(\gamma, m)=&&m(i(\gamma, 1)+p_-+p_0-r)+2\sum_{j=1}^r
E\left(\frac{m\theta_j}{2\pi}\right)-r
-p_--p_0\nn\\&&-\frac{1+(-1)^m}{2}(q_0+q_+)+2\left(
\sum_{j=1}^{r_*}\varphi\left(\frac{m\alpha_j}{2\pi}\right)-r_*\right).
\lb{3.12}\eea \bea \nu(\gamma, m)=&&\nu(\gamma,
1)+\frac{1+(-1)^m}{2}(q_-+2q_0+q_+)+2(r+r_*+r_0)\nn\\
&&-2\left(\sum_{j=1}^{r}\varphi\left(\frac{m\theta_j}{2\pi}\right)+
\sum_{j=1}^{r_*}\varphi\left(\frac{m\alpha_j}{2\pi}\right)
+\sum_{j=1}^{r_0}\varphi\left(\frac{m\beta_j}{2\pi}\right)\right).\lb{3.13}\eea
\bea \hat i(\gamma, 1)=i(\gamma, 1)+p_-+p_0-r+\sum_{j=1}^r
\frac{\theta_j}{\pi}.\lb{3.14}\eea We have $i(\gamma, 1)$ is odd if
$f(1)=N_1(1, 1)$, $I_2$, $N_1(-1, 1)$, $-I_2$, $N_1(-1, -1)$ and
$R(\theta)$; $i(\gamma, 1)$ is even if $f(1)=N_1(1, -1)$ and $
N_2(\omega, b)$; $i(\gamma, 1)$ can be any integer if $\sigma (f(1))
\cap \U=\emptyset$ by Theorems 8.1.4-8.1.7, 8.2.1, 8.2.3 and 8.2.4
of \cite{Lon3}.}

We have the following properties in the index iteration theory.

{\bf Theorem 3.4. } (cf. Theorem 2.4 of \cite{LoZ} ) {\it Let
$\gamma\in\P_\tau(2n)$ and $M = \gamma(\tau)$. Then for any
$m_1, m_2\in\N$, there holds
\bea&&\nu(\gamma, m_1)+\nu(\gamma, m_2)-\nu(\gamma, (m_1, m_2))-\frac{e(M)}{2}
\nn\\\le&&i(\gamma, m_1+m_2)-i(\gamma, m_1)-i(\gamma, m_2)
\nn\\\le &&\nu(\gamma, (m_1, m_2))-\nu(\gamma, m_1+m_2)+\frac{e(M)}{2},
\nn\eea
where $(m_1, m_2)$ is the greatest common divisor of $m_1$ and $m_2$.
}

The following is the enhanced common index jump theorem obtained in
\cite{DLW2}, which is a generalization of the standard common index
jump theorem in \cite{LoZ}.

{\bf Theorem 3.5.} (cf. Theorem 3.5 of \cite{DLW2}) {\it Let $\gamma_k\in\mathcal{P}_{\tau_k}(2n)$
for $k=1,\cdots,q$ be a finite collection of symplectic paths. Let
$M_k=\ga_k(\tau_k)$. We extend $\ga_k$ to $[0,+\infty)$ by
(\ref{3.9}) inductively. Suppose \bea  \hat{i}(\ga_k,1) > 0, \qquad
\forall\ k=1,\cdots,q.  \nn\eea We define \bea  \check{m} \equiv
\check{m}(\ga_1,\ldots,\ga_k) = \min\{\check{m}(M_k)\;|\;1\le k\le
q\}.  \nn\eea Then for every integer $\bar{m}\in \N$, there exist
infinitely many $(q+1)$-tuples $(N, m_1,\cdots,m_q) \in \N^{q+1}$
such that we have \bea  \nu(\ga_k,2m_k-m) = \nu(\ga_k,2m_k+m) =
\nu(\ga_k,1), \quad \forall\;1\le k\le q, 1\le m<\check{m},  \nn\eea
and the following hold for all $1\le k\le q$ and $1\le m\le
\bar{m}$, \bea
\nu(\ga_k,2m_k-m) &=& \nu(\ga_k,2m_k+m) = \nu(\ga_k, m),   \lb{3.15}\\
i(\ga_k,2m_k+m) &=& 2N+i(\ga_k,m),                         \lb{3.16}\\
i(\ga_k,2m_k-m) &=& 2N-i(\ga_k,m)-2(S^+_{M_k}(1)+Q_k(m)),  \lb{3.17}\\
i(\ga_k, 2m_k)&=& 2N -(S^+_{M_k}(1)+C(M_k)-2\Delta_k), \lb{3.18}\eea
where $C(M_k)=\sum_{0<\th<2\pi}S_{M_k}^-(e^{\sqrt{-1}\th})$ and \bea
\Delta_k &=&
\sum_{0<\{m_k\th/\pi\}<\delta}S^-_{M_k}(e^{\sqrt{-1}\th}),
\lb{3.19}\\ Q_k(m) &=&
\sum_{\theta\in(0,\,2\pi),\,e^{\sqrt{-1}\th}\in\sg(M_k),\atop
\{\frac{m_k\th}{\pi}\}
                   = \{\frac{m\th}{2\pi}\}=0}S^-_{M_k}(e^{\sqrt{-1}\th}). \lb{3.20}\eea
Moreover we have \bea
\min\left\{\left\{\frac{m_k\theta}{\pi}\right\},\,1-\left\{\frac{m_k\theta}{\pi}\right\}\right\}
<\delta,\lb{3.21}\eea whenever $e^{\sqrt{-1}\theta}\in\sigma(M_k)$
and $\delta$ can be chosen as small as we want (cf. (4.43) of
\cite{LoZ} ). More precisely, by (4.10) and (4.40) in \cite{LoZ}
 , we have
\bea m_k=\left(\left[\frac{N}{M\hat i(\gamma_k,
1)}\right]+\chi_k\right)M,\quad 1\le k\le q,\lb{3.22}\eea where
$\chi_k=0$ or $1$ for $1\le k\le q$ and $M\in\N$ such that $\frac{M\theta}{\pi}\in\Z$
whenever $e^{\sqrt{-1}\theta}\in\sigma(M_k)$ and
$\frac{\theta}{\pi}\in\Q$ for some $1\le k\le q$. Furthermore, given
$M_0\in\N$, by the proof of Theorem 4.1 of \cite{LoZ}, we may
further require  $M_0|N$ (since the closure of the set $\{\{Nv\}:
N\in\N, \;M_0|N\}$ is still a closed additive subgroup of $\bf T^h$
for some $h\in\N$, where we use notations  as (4.21)  in \cite{LoZ}.
  Then we can use the proof of Step 2 in Theorem 4.1 of \cite{LoZ}
 to get $N$).

 In fact,  Let  $\mu_i=\sum_{\theta\in(0, 2\pi)}S_{M_i}^-(e^{\sqrt{-1}\theta})$ for $1\le i\le q$
and $\alpha_{i, j}=\frac{\theta_j}{\pi}$ where
$e^{\sqrt{-1}\theta_j}\in\sigma(M_i)$ for $1\le j\le \mu_i$ and
$1\le i\le q$.  Let $h=q+\sum_{1\le i\le q}\mu_i$ and \bea
v=\left(\frac{1}{M\hat i(\gamma_1, 1)},\dots, \frac{1}{M\hat
i(\gamma_q, 1)}, \frac{\alpha_{1, 1}}{\hat i(\gamma_1, 1)},
\frac{\alpha_{1, 2}}{\hat i(\gamma_1, 1)}, \dots\frac{\alpha_{1,
\mu_1}}{\hat i(\gamma_1, 1)}, \frac{\alpha_{2, 1}}{\hat i(\gamma_2,
1)},\dots, \frac{\alpha_{q, \mu_q}}{\hat i(\gamma_q,
1)}\right)\in\R^h. \lb{3.23}\eea Then  the above theorem is
equivalent to find a vertex \bea\chi=(\chi_1,\ldots,\chi_q,\chi_{1,
1},\chi_{1, 2},\ldots,\chi_{1, \mu_1}, \chi_{2, 1},\ldots,\chi_{q,
\mu_q})\lb{3.24}\eea of the cube $[0,\,1]^h$ and infinitely many
integers $N\in\N$ such that \bea|\{Nv\}-\chi|<\epsilon\lb{3.25}\eea
 for any given $\epsilon$ small enough (cf. P. 346 and 349 of \cite{LoZ}
).}
 \medskip

The following theorem implies that  we
can choose two special tuples of
 integers $(N, m_1, \ldots,m_q)$ and
$(N^\prime, m_1^\prime,\cdots,m_q^\prime)\in \N^{q+1}$ in Theorem
3.5 such that the indices have some symmetric property.

\medskip

{\bf Theorem 3.6.} (cf. Theorem 2.8 of \cite{HW}) {\it Let $\gamma_k\in\mathcal{P}_{\tau_k}(2n)$
for $k=1,\cdots,q$ be a finite collection of symplectic paths. Let
$M_k=\ga_k(\tau_k)$. We extend $\ga_k$ to $[0,+\infty)$ by
(\ref{3.9}) inductively. Suppose \be  \hat{i}(\ga_k,1) > 0, \qquad
\forall\ k=1,\cdots,q.  \lb{3.26}\ee Then for every integer
$\bar{m}\in \N$, there exist infinitely many $(q+1)$-tuple pairs
$(N, m_1,\cdots,m_q)\in\N^{q+1}$ and  $(N^\prime,
m_1^\prime,\cdots,m_q^\prime)\in \N^{q+1}$ such that the following
hold for all $1\le k\le q$ and $1\le m\le \bar{m}$, \bea
\nu(\ga_k,2m_k-m) &=& \nu(\ga_k,2m_k+m) = \nu(\ga_k, m),   \lb{3.27}\\
i(\ga_k,2m_k+m) &=& 2N+i(\ga_k,m),                         \lb{3.28}\\
i(\ga_k,2m_k-m) &=& 2N-i(\ga_k,m)-2(S^+_{M_k}(1)+Q_k(m)),  \lb{3.29}\\
i(\ga_k, 2m_k)&=& 2N -(S^+_{M_k}(1)+C(M_k)-2\Delta_k),     \lb{3.30}\\
\nu(\ga_k,2m_k^\prime-m) &=& \nu(\ga_k,2m_k^\prime+m) = \nu(\ga_k, m),   \lb{3.31}\\
i(\ga_k,2m_k^\prime+m) &=& 2N^\prime+i(\ga_k,m),                         \lb{3.32}\\
i(\ga_k,2m_k^\prime-m) &=& 2N^\prime-i(\ga_k,m)-2(S^+_{M_k}(1)+Q_k^\prime(m)),  \lb{3.33}\\
i(\ga_k, 2m_k^\prime)&=& 2N^\prime
-(S^+_{M_k}(1)+C(M_k)-2\Delta_k^\prime),     \lb{3.34}\eea where
$C(M_k)=\sum_{0<\th<2\pi}S_{M_k}^-(e^{\sqrt{-1}\th})$ and \bea
\Delta_k &=& \sum_{0<\{m_k\th/\pi\}<\delta}S^-_{M_k}(e^{\sqrt{-1}\th}), \lb{3.35}\\
\Delta_k^\prime &=& \sum_{0<\{m_k^\prime\th/\pi\}<\delta}S^-_{M_k}(e^{\sqrt{-1}\th}), \lb{3.36}\\
 Q_k(m)& = &\sum_{\theta\in(0,\,2\pi),\,e^{\sqrt{-1}\th}\in\sg(M_k),\atop \{\frac{m_k\th}{\pi}\}
                   = \{\frac{m\th}{2\pi}\}=0}S^-_{M_k}(e^{\sqrt{-1}\th}),\lb{3.37}\\
 Q_k^\prime(m)& = &\sum_{\theta\in(0,\,2\pi),\,e^{\sqrt{-1}\th}\in\sg(M_k),\atop \{\frac{m_k^\prime\th}{\pi}\}
                   = \{\frac{m\th}{2\pi}\}=0}S^-_{M_k}(e^{\sqrt{-1}\th}). \lb{3.38}\eea
 Furthermore, $\Delta_k$  and $\Delta_k'$  satisfy the following
relationship \be \Delta_k + \Delta_k^\prime
=\sum_{\frac{\th}{2\pi}\in (0, \,1)\setminus
\Q}S_{M_k}^-(e^{\sqrt{-1}\th}) ,\qquad\forall\ 1\le k\le q.
\lb{3.39}\ee }

\medskip

\setcounter{figure}{0}
\setcounter{equation}{0}
\section{Proof of Theorem 1.2}

In this section, we give the proof of Theorem 1.2.

Firstly we have the following index estimates for closed geodesics.

{\bf Theorem 4.1.} (cf. Theorem 1.2 of \cite{DLW2}) {\it Let $(M,F)$ be an even dimensional compact simply-connected manifold with a
bumpy irreversible Finsler metric $F$. Then every prime closed geodesic $c$ on $(M,F)$ satisfies
$i(c)>0$, provided the flag curvature $K$ of the Chern connection satisfies $K\ge 0$.}

{\bf Lemma 4.2.} {\it Let $c$ be a closed geodesic on a
compact simply-connected Finsler manifold $(M,F)$
satisfying $\hat i(c)>0$. Then there exists $\tilde M\in\N$
such that
\bea i(c^{m+l})\ge i(c^l)+\dim M,\quad \forall l\in\N,\;m\ge \tilde m.\lb{4.1}\eea
}

{\bf Proof.}
Firstly by Theorem 1.4 of \cite{Rad1},
we have
\bea |i(c^m)-m\hat i(c)|\le \dim M-1,\quad m\in\N.\lb{4.2}\eea
By Theorem 3.4,
\bea i(c^{m+l})-i(c^l)&\ge& i(c^m)+\nu(c^m)+\nu(c^l)
-\nu(c^{(m,\, l)})-\frac{e(P_c)}{2}
\nn\\&\ge&i(c^m)+0+0-(2\dim M-2)-(\dim M-1)=i(c^m)-3\dim M+3
\nn\\&\ge&m\hat i(c)-(\dim M-1)-3\dim M+3=m\hat i(c)-4\dim M+4.\lb{4.3}\eea
where the second inequality follows from $0\le \nu(d)\le 2\dim M-2$
and $0\le e(P_d)\le 2\dim M-2$ for any closed geodesic $d$;
the last inequality follows from (\ref{4.2}).
Since $\hat i(c)>0$, we can choose $\tilde m\in\N$ satisfies
$\tilde m\hat i(c)-4\dim M+4>\dim M$,
then the lemma holds.\hfill\hb

When the Finsler metric $F$ is bumpy, for every prime closed geodesic $c$, by Theorem 3.2 the
basic normal form decomposition of the linearized Poincar\'{e} map $P_{c}$  possesses the following
form
\bea f_c(1)
&=& R(\th_1)\,\dm\,\cdots\,\dm\,R(\th_r)\,\dm\,D(\lm_1)\dm \,\cdots\,\dm D(\lm_s) \nn\\
& & \dm\,N_2(e^{\sqrt{-1}\aa_{1}},u_{1})\,\dm\,\cdots\,\dm\,N_2(e^{\sqrt{-1}\aa_{r_{\ast}}},u_{r_{\ast}})\,
\nn\\ & &   \dm\,N_2(e^{\sqrt{-1}\bb_{1}},v_{1})\,\dm\,\cdots\,\dm\,N_2(e^{\sqrt{-1}\bb_{r_{0}}},v_{r_{0}}), \lb{4.4}\eea
where $N_2(e^{\sqrt{-1}\aa_{j}}, u_j)$'s
are non-trivial and   $N_2(e^{\sqrt{-1}\beta_{j}}, v_j)$'s  are trivial basic normal
forms; $\lm_j\in \{\pm 2\}$ for $1\le j\le s$, $\frac{\th_{j}}{2\pi}\not\in\Q$ for $1\le j\le r$,
$\frac{\aa_{j}}{2\pi}\not\in\Q$ for $1\le j\le r_{\ast}$, $\frac{\bb_{j}}{2\pi}\not\in\Q$ for $1\le j\le r_0$,
and
\be r + s + 2r_{\ast} + 2r_0 = \dim M - 1. \lb{4.5}\ee
Then by Theorem 3.3, we obtain the index iteration formula of $c^m$
\be i(c^m) = m(i(c)-r) + 2\sum_{j=1}^r\left[\frac{m\th_j}{2\pi}\right] +r,\quad \nu(c^m)=0,
    \qquad\forall\ m\ge 1, \lb{4.6}\ee
where we have used the fact that $E(x)=[x]+1$ for $x\notin\Z$.

{\bf Proof of Theorem 1.2.}
Let $(M,F)$ be a compact simply-connected manifold with a bumpy
irreversible Finsler metric $F$ and satisfy
$H^*(M;\Q)\cong T_{d,n+1}(x)$ for some integers
$d\in 2\N$ and  $n\in\N$. Suppose the flag curvature $K\ge 0$,
then by Theorem 4.1, we have $i(c)>0$ for any closed geodesic $c$
on $(M,F)$.

We make the following assumption.

\medskip

{\bf (FCG)} {\it Suppose that there exist only finitely many prime closed geodesics
$\{c_k\}_{1\le k\le q}$ with $i(c_k)\ge 1$ for $k=1,\cdots,q$ on $(M,F)$.}

\medskip

We carry out the proof in six steps.

{\bf Step 1}. {\it There exist $\frac{dn(n+1)}{4}$
closed geodesics on $(M, F)$.}

By the assumption (FCG), there exist only finitely many distinct prime closed geodesics
on  $(M,F)$, then any closed geodesic $c_k$ among $\{c_k\}_{k=1}^q$ must have
positive mean index (cf. Theorem 3 of \cite{BaK} or Lemma 3.2 of \cite{DuL2}), i.e.,
\be \hat{i}(c_k)>0,\qquad 1\le k\le q. \lb{4.7}\ee
In fact, by Corollary 1.3 of \cite{Rad1}, $\hat i(c_k)=0$
implies $i(c_k^m)=0$ for $m\in\N$, contradict to $i(c_k)>0$,
this gives another proof of (\ref{4.7}).
Thus by Lemma 4.2, we have $\tilde m_k\in\N$ for $1\le k\le q$ such that
\bea i(c_k^{m+l})\ge i(c_k^l)+\dim M,\quad \forall l\in\N,\;m\ge \tilde m_k.
\lb{4.8}\eea
Denote by $\bar{m}=\max_{1\le k\le q}\{\tilde m_k\}+2$.

By Theorem 3.5, there exist infinitely many $(q+1)$-tuples $(N, m_1, \cdots, m_q)\in\N^{q+1}$ such that
$\bar{m}+2\le\min_{1\le k\le q}\{2m_k\}$ and
for any $1\le k\le q$, there holds
\bea
i(c_k^{2m_k-m}) &=& 2N-i(c_k^m),\quad 1\le m\le\bar{m}, \lb{4.9}\\
i(c_k^{2m_k}) &=& 2N-C(P_{c_k})+2\Delta_k,  \lb{4.10}\\
i(c_k^{2m_k+m}) &=& 2N+i(c_k^m),\quad 1\le m\le\bar{m}, \lb{4.11}\eea
where (\ref{4.9}) and (\ref{4.10}) follows from
(\ref{3.17}) and (\ref{3.18}), $S^+_{P_{c_k}}(1)=0$,
$Q_k(m)=0$ for $m\in\N$ since $F$ is bumpy.

By (\ref{4.8}), we have
\bea i(c_k^{2m_k+m})&=&i(c_k^{(2m_k+1)+(m-1)})\ge i(c_k^{2m_k+1})+\dim M,
\quad \forall m\ge \bar m+1,\lb{4.12}\\
i(c_k^{2m_k-m})&\le&i(c_k^{(2m_k-m)+(m-1)})= i(c_k^{2m_k-1}),
\quad  \bar m+1\le m\le 2m_k-1,\lb{4.13}\eea
On one hand, there holds $i(c_k^m)\ge i(c_k)$ for any $m\ge 1$ by the Bott-type formulae (cf.
\cite{Bot} and Theorem 9.2.1 of \cite{Lon3}).
Thus by (\ref{4.9})-(\ref{4.13}) and the
assumption $i(c_k)\ge 1$ for $1\le k\le q$ in (FCG),
we have
\bea
&& i(c_k^{2m_k-m})\le 2N-i(c_k)\le 2N-1,\qquad \forall\ 1\le m< 2m_k, \lb{4.14}\\
&& i(c_k^{2m_k})=2N-C(P_{c_k})+2\Delta_k,\lb{4.15}\\
&& i(c_k^{2m_k+m})\ge 2N+i(c_k)\ge 2N+1,\qquad \forall\ m\ge 1.\lb{4.16}\eea

{\bf Claim 1.} {\it For $N\in\N$ in Theorem 3.5 satisfying
(\ref{4.14})-(\ref{4.16}) and $2NB(d,n)\in2\Z$,
we have
\be \sum_{1\le k\le q} 2m_k\gamma_{c_k}=2NB(d,n).\lb{4.17}\ee}

In fact, here we follow the ideas in \cite{Wan1} and \cite{Wan2} to prove this Claim. By Lemma 2.3,
(\ref{3.22}) and (\ref{3.25}), it yields
\bea \left|2NB(d,n)-\sum_{k=1}^q 2m_k\ga_{c_k}\right|
&=& \left|\sum_{k=1}^q\frac{2N\ga_{c_k}}{\hat{i}(c_k)}-\sum_{k=1}^q 2\ga_{c_k}
                    \left(\left[\frac{N}{M\hat{i}(c_k)}\right]+\chi_k\right)M\right| \nn\\
&=& \left|\sum_{k=1}^q\frac{2N\ga_{c_k}}{\hat{i}(c_k)}-\sum_{k=1}^q 2\ga_{c_k}
                    \left(\frac{N}{M\hat{i}(c_k)}-\left\{\frac{N}{M\hat{i}(c_k)}\right\}+\chi_k\right)M\right| \nn\\
&\le& 2M\sum_{k=1}^q |\ga_{c_k}|\left|\left\{\frac{N}{M\hat{i}(c_k)}\right\}-\chi_k\right| \nn\\
&<& 2M\ep\sum_{k=1}^q|\ga_{c_k}|< 1. \lb{4.18}\eea
provided  $\ep<\frac{1}{1+2M\sum_{1\le k\le q}|\ga_{c_k}|}$.
Since $2NB(d, n)$ and  $2m_k\ga_{c_k}$'s are integers, Claim 1 is proved.

\medskip

Now by Lemma 2.1, Definition 2.2 and Theorem 3.3, it yields
\bea \sum_{m=1}^{2m_k} (-1)^{i(c_k^m)} \dim \ol{C}_{i(c_k^{m})}(E,c_k^m)
&=& \sum_{i=0}^{m_k-1} \sum_{m=2i+1}^{2i+2} (-1)^{i(c_k^m)} \dim \ol{C}_{i(c_k^{m})}(E,c_k^m) \nn\\
&=& \sum_{i=0}^{m_k-1} \sum_{m=1}^{2} (-1)^{i(c_k^m)} \dim \ol{C}_{i(c_k^{m})}(E,c_k^m) \nn\\
&=& m_k \sum_{m=1}^{2} (-1)^{i(c_k^m)} \dim \ol{C}_{i(c_k^{m})}(E,c_k^m) \nn\\
&=& 2m_k\ga_{c_k},\qquad \forall\ 1\le k\le q, \lb{4.19}\eea
where the second equality follows from Lemma 2.1 and the fact
$i(c_k^{m+2})-i(c_k^m)\in 2\Z$ for
$m\in\N$ from Theorem 3.3, and the last equality follows from Lemma 2.1 and Definition 2.2.

By (\ref{4.16}) and Lemma 2.1, we know that all $c_k^{2m_k+m}$'s with $m\ge1$ and $1\le k\le q$
have no contributions to the alternative sum $\sum_{p=0}^{2N}(-1)^p M_p$. Similarly again by Lemma 2.1
and (\ref{4.14}), all $c_k^{2m_k-m}$'s with $1\le m<2m_k$ and $1\le k\le q$ only have contributions
to $\sum_{p=0}^{2N}(-1)^p M_p$.

Thus for the Morse-type numbers $M_p$'s defined by (\ref{2.9}),
by (\ref{4.19}) we have
\bea \sum_{p=0}^{2N}(-1)^p M_p
&=& \sum_{k=1}^{q}\ \sum_{1\le m\le 2m_k, \atop i(c_k^{m})\le 2N} (-1)^{i(c_k^m)} \dim \ol{C}_{i(c_k^{m})}(E,c_k^m) \nn\\
&=& \sum_{k=1}^{q}\ \sum_{m=1}^{2m_k} (-1)^{i(c_k^m)}\dim\ol{C}_{i(c_k^{m})}(E,c_k^m) \nn\\
& & -\sum_{1\le k\le q, \atop i(c_k^{2m_k})\ge 2N+1} (-1)^{i(c_k^{2m_k})} \dim \ol{C}_{i(c_k^{2m_k})}(E,c_k^{2m_k}) \nn\\
&=& \sum_{k=1}^{q} 2m_k\ga_{c_k}
  -\sum_{1\le k\le q, \atop i(c_k^{2m_k})\ge 2N+1} (-1)^{i(c_k^{2m_k})} \dim \ol{C}_{i(c_k^{2m_k})}(E,c_k^{2m_k}).
      \lb{4.20}\eea

Denote by
\bea
\Omega_+^e &=& \{1\le k\le q\ |\ i(c_k^{2m_k})\ge 2N+1,\ i(c_k^{2m_k})-i(c_k)\in 2\N_0,\ i(c_k)\in 2\N\}, \lb{4.21}\\
\Omega_+^o &=& \{1\le k\le q\ |\ i(c_k^{2m_k})\ge 2N+1,\ i(c_k^{2m_k})-i(c_k)\in 2\N_0,\ i(c_k)\in 2\N-1\}, \lb{4.22}\\
\Omega_-^e &=& \{1\le k\le q\ |\ i(c_k^{2m_k})\le 2N-1,\ i(c_k^{2m_k})-i(c_k)\in 2\N_0,\ i(c_k)\in 2\N\}, \lb{4.23}\\
\Omega_-^o &=& \{1\le k\le q\ |\ i(c_k^{2m_k})\le 2N-1,\ i(c_k^{2m_k})-i(c_k)\in 2\N_0,\ i(c_k)\in 2\N-1\}, \lb{4.24}\eea
and $^\#\Omega_\pm^e=N_\pm^e$, $^{\#}\Omega_\pm^o=N_\pm^o$.
Note that by Theorem 3.5, we can suppose that $N$ is large enough and is a multiple of $D=d(n+1)-2$.

Thus we have
\bea -\frac{Ndn(n+1)}{D}+N_+^o-N_+^e
&=& 2NB(d,n)+N_+^o-N_+^e  \nn\\
&=& \sum_{k=1}^q 2m_k\gamma_{c_k}+N_+^o-N_+^e \nn\\
&=& \sum_{p=0}^{2N}(-1)^p M_p \nn\\
&\ge&\sum_{p=0}^{2N}(-1)^p b_p(\Lm M).
\lb{4.25}\eea
In fact, the first equality follows from Theorem 2.3;
the second equality follows from Claim 1;
the third equality follows from (\ref{4.20})-(\ref{4.22})
and Lemma 2.1; the last inequality follows from Theorem 2.5.

Write $N=DT$ for some $T\in\N$,
then by Lemma 2.4,
\bea&&\sum_{p=0}^{2N}(-1)^p b_p(\Lm M) \nn\\
&=&\sum_{0\le p\le 2N,\atop{p\in2\Z}}(-1)^p b_p(\Lm M)
+\sum_{0\le p\le d-2,\atop{p\in2\N-1}}(-1)^p b_p(\Lm M)
\nn\\&&+\sum_{d-1\le p< nd-1,\atop{p\in2\N-1}}(-1)^p b_p(\Lm M)
+\sum_{nd-1\le p\le2N,\atop{p\in\Omega(d,n)}}(-1)^p b_p(\Lm M)
+\sum_{nd-1\le p\le2N,\atop{p\in2\N-1, p\notin\Omega(d,n)}}(-1)^p b_p(\Lm M)
\nn\\
&=& 0+0+S_1+S_2+S_3.\lb{4.26}\eea
Write $p=d-1+kd+2s$ for $k,s\in\Z$ with $0\le k\le n-2$
and $0\le 2s\le d-2$, then
\bea S_1&=&\sum_{d-1\le p< nd-1,\atop{p\in2\N-1}}(-1)^p b_p(\Lm M)
=\sum_{0\le k\le n-2, \atop{0\le 2s\le d-2}}(-1)^{d-1+kd+2s} b_{d-1+kd+2s}(\Lm M)
\nn\\&=&\sum_{0\le k\le n-2, \atop{0\le 2s\le d-2}}
(-1)\left(\left[\frac{d-1+kd+2s-(d-1)}{d}\right]+1\right)
\nn\\&=&-\sum_{0\le k\le n-2, \atop{0\le 2s\le d-2}}
(k+1)=-\frac{n(n-1)d}{4}.
\lb{4.27}\eea
Write $p=d-1+iD+jd$ for $p\in\Omega(d, n)$ and $N=DT$ for some
$T\in\N$, then
\bea
S_2&=&\sum_{nd-1\le p\le2N,\atop{p\in\Omega(d,n)}}(-1)^p b_p(\Lm M)
=\sum_{nd-1\le d-1+iD+jd\le2DT,\atop{0\le j\le n-1}}(-1)^p b_p(\Lm M)
\nn\\&=&\sum_{1\le i\le2T-1,\atop{0\le j\le n-1}}(-(n+1))
=-(2T-1)n(n+1),\lb{4.28}\eea
where we use the fact that $D=(n+1)d-2$,
then $nd-1\le d-1+D=d-1+(n+1)d-2=(n+2)d-3$ and
$d-1+(2T-1)D+(n-1)d=2TD+nd-1-D=2TD+nd-1-((n+1)d-2)=
2TD-d+1\le 2TD$ since $d\ge 2$, hence the ranges of $i, j$ are $1\le i\le 2T-1$
and $0\le j\le n-1$.
\bea
S_3&=&\sum_{nd-1\le p\le2N,\atop{p\in2\N-1, p\notin\Omega(d,n)}}(-1)^p b_p(\Lm M)
\sum_{nd-1\le p\le2DT,\atop{p\in2\N-1, p\notin\Omega(d,n)}}(-n)
\nn\\&=&-n\left(\frac{2DT-1-(nd-1)}{2}+1-(2T-1)n\right)
\nn\\&=&-n\left(DT-\frac{nd}{2}-(2T-1)n+1\right)
\nn\\&=&-n\left(((n+1)d-2)T-\frac{nd}{2}-(2T-1)n+1\right),
\lb{4.29}\eea
where the second equality follows since the number
$\{p| nd-1\le p\le2N,p\in\Omega(d,n)\}$ is $(2T-1)n$
by the computation of $S_2$.

Therefore
\bea&&\sum_{p=0}^{2N}(-1)^p b_p(\Lm M)=S_1+S_2+S_3
\nn\\=&&-\frac{n(n-1)d}{4}-(2T-1)n(n+1)-n\left((n+1)d-2)T-\frac{nd}{2}-(2T-1)n+1\right)
\nn\\=&&-n(n+1)dT+\frac{n(n+1)d}{4}
\nn\\=&&2NB(d, n)+\frac{n(n+1)d}{4},\lb{4.30}\eea
where the last equality holds by (\ref{2.4}).

From (\ref{4.25}) and (\ref{4.30}), we obtain
\be N_+^o\ge \frac{dn(n+1)}{4}, \lb{4.31}\ee
i.e., we find  $\frac{dn(n+1)}{4}$ closed geodesics $c_{k_i}$ with
$\{k_1, k_2,\ldots, k_\frac{dn(n+1)}{4}\}\subset\Omega^o_+$.
This proves Step 1.
\medskip

{\bf Step 2}. {\it There exist another $\frac{dn(n+1)}{4}$ closed
geodesics on $(M, F)$ which are distinct from those in Step 1.}

We use Theorem 3.6 to get infinitely many
 $(q+1)$-tuples $(N, m_1, \cdots, m_q),\,(N', m_1', \cdots, m_q')\in\N^{q+1}$
such that (\ref{4.14})-(\ref{4.16}) hold and for every
$1\le k\le q$, there hold
\bea
&& i(c_k^{2m_k'-m})\le 2N'-1,\qquad \forall\ 1\le m<2m'_k, \lb{4.32}\\
&& i(c_k^{2m_k'})=2N'-C(P_{c_k})+2\Delta_k',  \lb{4.33}\\
&& i(c_k^{2m_k'+m})\ge 2N'+1,\qquad \forall\ m\ge 1,  \lb{4.34}\eea
Moreover, \bea
\Delta_k &=& \sum_{0<\{m_k\th/\pi\}<\delta}S^-_{P_{c_k}}(e^{\sqrt{-1}\th}), \lb{4.35}\\
\Delta_k^\prime &=&
\sum_{0<\{m_k^\prime\th/\pi\}<\delta}S^-_{P_{c_k}}(e^{\sqrt{-1}\th}),
\lb{4.36} \eea satisfies \be \Delta_k + \Delta_k^\prime
=\sum_{\frac{\th}{2\pi}\in (0, \,1)\setminus
\Q}S_{P_{c_k}}^-(e^{\sqrt{-1}\th})
=\sum_{0<\th<2\pi}S_{P_{c_k}}^-(e^{\sqrt{-1}\th})=C(P_{c_k}),
\quad\forall\ 1\le k\le q, \lb{4.37}\ee where the second equality
holds since $e^{\sqrt{-1}\th}\notin\sigma(P_c)$ for
$\frac{\th}{2\pi}\in (0, \,1)\cap\Q$ by (\ref{4.4}).

Similarly to Step 1, we define
\bea
\Omega_+^{'e} &=& \{1\le k\le q\ |\ i(c_k^{2m_k'})\ge 2N'+1,\ i(c_k^{2m_k'})-i(c_k)\in 2\N_0,\ i(c_k)\in 2\N\},   \lb{4.38}\\
\Omega_+^{'o} &=& \{1\le k\le q\ |\ i(c_k^{2m_k'})\ge 2N'+1,\ i(c_k^{2m_k'})-i(c_k)\in 2\N_0,\ i(c_k)\in 2\N-1\}, \lb{4.39}\\
\Omega_-^{'e} &=& \{1\le k\le q\ |\ i(c_k^{2m_k'})\le 2N'-1,\ i(c_k^{2m_k'})-i(c_k)\in 2\N_0,\ i(c_k)\in 2\N\},   \lb{4.40}\\
\Omega_-^{'o} &=& \{1\le k\le q\ |\ i(c_k^{2m_k'})\le 2N'-1,\ i(c_k^{2m_k'})-i(c_k)\in 2\N_0,\ i(c_k)\in 2\N-1\}, \lb{4.41}\eea
and $^{\#}\Omega_\pm^{'e}=N_\pm^{'e}$, $^{\#}\Omega_\pm^{'o}=N_\pm^{'o}$.

By (\ref{4.33}) and (\ref{4.37}) it yields
\be i(c_k^{2m_k'}) = 2N'-C(P_{c_k})+2(C(P_{c_k})-\Delta_k)
=2N'+C(P_{c_k})-2\Delta_k. \lb{4.42}\ee

By (\ref{4.15}) and (\ref{4.42}), we have \bea \Omega_{\pm}^e =
\Omega_{\mp}^{'e}, \qquad \Omega_{\pm}^o = \Omega_{\mp}^{'o},
\lb{4.43}\\
N_{\pm}^e = N_{\mp}^{'e}, \qquad N_{\pm}^o = N_{\mp}^{'o}.
\lb{4.44}\eea

As in Step 1, we have
\be  N_+^{'o}\ge\frac{dn(n+1)}{4}, \lb{4.45}\ee
i.e., we find another $\frac{dn(n+1)}{4}$ closed geodesics $c_{k_i}'$ with
$\{k_1', k_2',\ldots, k_\frac{dn(n+1)}{4}'\}
\subset\Omega^{'o}_+=\Omega^o_-$.
By (\ref{4.22}) and (\ref{4.24}), $\Omega^{o}_+\cap\Omega^o_-=\emptyset$,
thus these closed geodesics are distinct from those obtained in Step 1.
So by (\ref{4.31}), (\ref{4.44}) and (\ref{4.45}) we get
\be q \ge N_+^{o}+N_-^{o} \ge \frac{dn(n+1)}{4}+\frac{dn(n+1)}{4}
= \frac{dn(n+1)}{2}. \lb{4.46}\ee

\medskip
{\bf Step 3.} {\it Suppose there are exactly $\frac{dn(n+1)}{2}=q$
closed geodesics, then we have $M_p=b_p(\Lambda M)$ for $p\in\Z$.}

By the assumption $q=\frac{dn(n+1)}{2}$ and (\ref{4.46}),
we have $N_+^{o}=N_-^{o} = \frac{dn(n+1)}{4}$,
and then
\be\Omega_+^{o}\cup\Omega_-^{o}=\{1, 2,\ldots,q\}.\lb{4.47}\ee
Hence by the definitions of $\Omega_\pm^{o}$, cf. (\ref{4.22}) and (\ref{4.24}),
\bea i(c_k)\in 2\N-1,\quad i(c_k^{2m_k})-i(c_k)\in 2\N_0,\qquad1\le k\le q.
\lb{4.48}\eea

By Theorem 3.3,
\bea i(c^{m+2})-i(c^m)= 0\;\mod\; 2,\quad m\in\N,\lb{4.49}\eea
for any closed geodesic $c$ on $(M, F)$.
Hence we have
\bea i(c_k^m)\in 2\N-1, \quad 1\le k\le q,\;\forall m\in\N.
\lb{4.50}\eea
Thus we have $M_{2k}=0=b_{2k}(\Lambda M)$ for $k\in\Z$ by Lemmas 2.1 and 2.4,
and then we have $M_{2k-1}=b_{2k-1}(\Lambda M)$ for $k\in\Z$ by Theorem 2.5.
This proves step 3.

{\bf Step 4.} {\it There exists a closed geodesic $c_{k_\ast}$
satisfying $-C(P_{c_{k_\ast}})+2\Delta_{k_\ast}=nd-1$ among the
$\frac{dn(n+1)}{4}$ closed geodesics found in Step 1.}

By a permutation of $\{1, 2,\ldots, q\}$, we may denote these
$\frac{dn(n+1)}{4}=\frac{q}{2}$ closed geodesics by $c_1, c_2,
\ldots c_{\frac{q}{2}}$, i.e., $\Omega^o_+=\{1,
2,\ldots,\frac{q}{2}\}$, and then $\Omega^o_-=\{\frac{q}{2}+1,
\frac{q}{2}+2,\ldots,q\}$ by (\ref{4.47}).

We prove Step 4 by contraction, i.e., assume
$-C(P_{c_k})+2\Delta_k<nd-1$ for $1\le k\le \frac{q}{2}$.
Then for $1\le k\le q$, we have
\bea
&& i(c_k^{2m_k-m})\le 2N-1,\qquad \forall\ 1\le m< 2m_k, \lb{4.51}\\
&& i(c_k^{2m_k})=2N-C(P_{c_k})+2\Delta_k<2N+(nd-1),\lb{4.52}\\
&& i(c_k^{2m_k+m})=2N+i(c_k^m),\qquad 1\le m\le\bar{m},\lb{4.53}\\
&& i(c_k^{2m_k+m})\ge 2N+1+\dim M,\qquad \forall m>\bar{m},\lb{4.54}\eea

In fact, (\ref{4.51}) holds by (\ref{4.14});
(\ref{4.52}) holds for $1\le k\le \frac{q}{2}$
by the assumption $-C(P_{c_k})+2\Delta_k<nd-1$
for $1\le k\le \frac{q}{2}$,
(\ref{4.52}) holds for $\frac{q}{2}+1\le k\le q$
by the definition of $\Omega^o_-$ in (\ref{4.24});
(\ref{4.53}) holds by (\ref{4.11});
(\ref{4.54}) follows from $i(c_k^{2m_k+m})\ge i(c_k^{2m_k+1})+\dim M
=2N+i(c_k)+\dim M\ge 2N+1+\dim M$ by (\ref{4.11}) and (\ref{4.12}).

By Step 3, we have \bea &&M_{2N+1}+M_{2N+3}+\cdots+M_{2N+nd-3}
\nn\\=&&
b_{2N+1}(\Lambda M)+b_{2N+3}(\Lambda M)+\cdots+b_{2N+nd-3}(\Lambda M).
\lb{4.55} \eea

Write $N=DT$, $p=2N+1+kd+2s$ for $k,s\in\Z$ with $0\le k\le n-1$
and $0\le 2s\le d-2$, we have
 \bea &&b_{2N+1}(\Lambda M)+b_{2N+3}(\Lambda M)+\cdots+b_{2N+nd-3}(\Lambda M)\nn\\
=&&\sum_{2N+1\le p\le2N+nd-3} b_p(\Lm M)
 =\sum_{2DT+1\le p\le2DT+nd-3} b_p(\Lm M)
\nn\\=&&\sum_{0\le k\le n-1}\sum_{0\le 2s\le d-2}
b_{2DT+1+kd+2s}(\Lm M)- b_{2DT+nd-1}(\Lm M)
\lb{4.56}\\=&&\sum_{0\le k\le n-1}\left(n\frac{d-2}{2}+(n+1)\right)-(n+1)
=\frac{dn^2}{2}-1.\lb{4.57}\eea
where the fourth equality holds from
$2DT+1+kd+d-2=2DT+d-1+kd\in\Omega(d, n)$;
$2DT+1>d-1+(2T-1)D+jd$ and
$2DT+nd-1<d-1+(2T+1)D+jd$ for $0\le j\le n-1$ since $D=(n+1)d-2$ and $d\ge 2$;
hence in each inner summation of (\ref{4.56}), there is
exactly one term belongs to $\Omega(d, n)$, the other $\frac{d-2}{2}$ terms
 do not belong to $\Omega(d, n)$.

By Lemma 2.4 and (\ref{4.27}), we have
\bea \sum_{0\le p< nd-1,\atop{p\in2\N-1}} b_p(\Lm M)
=\frac{n(n-1)d}{4}.\lb{4.58}
\eea
Therefore by Step 3,
\bea M_{1}+M_{3}+\cdots+M_{nd-3}=
b_{1}(\Lambda M)+b_{3}(\Lambda M)+\cdots+b_{nd-3}(\Lambda M)=\frac{n(n-1)d}{4}.\lb{4.59} \eea
Hence by Lemma 2.1, there are
pairwise distinct $\{(k_l, m_l)\}_{1\le l\le\frac{n(n-1)d}{4}}=\Xi$
such that
\bea i(c_k^m)\in \left\{\matrix{
     \{d-1, d+1, \ldots, nd-3\},&\quad {\it if}\;\; (k, m)\in\Xi, \cr
     \{nd-1, nd+1,\cdots\}, &\quad {\it otherwise}, \cr}\right.
\lb{4.60}\eea
where $i(c_k^m)\ge d-1$ follows from $M_k=b_k(\Lambda M)=0$
for $k<d-1$ by Step 3 and Lemma 2.4.

By (\ref{4.53}) and (\ref{4.54}),
\bea i(c_{k_l}^{2m_{k_l}+m_l})=2N+i(c_{k_l}^{m_l}),\qquad
1\le l\le \frac{n(n-1)d}{4}.\lb{4.61}\eea
In fact, $m_l\le \bar m$ must hold, since otherwise
$i(c_{k_l}^{m_l})\ge i(c_{k_l})+\dim M\ge 1+\dim M=1+dn$ by (\ref{4.8}),
contradicts to $i(c_{k_l}^{m_l})\le dn-3$. Hence we have
\bea i(c_{k_l}^{2m_{k_l}+m_l})\in\{2N+d-1, 2N+d+1,\ldots, 2N+nd-3\},
\qquad
1\le l\le \frac{n(n-1)d}{4}.\lb{4.62}\eea

Hence the closed geodesics $c_{k_l}^{2m_{k_l}+m_l}$'s for $1\le l\le \frac{n(n-1)d}{4}$
contribute $\frac{n(n-1)d}{4}$ to the summation
$M_{2N+1}+M_{2N+3}+\cdots+M_{2N+nd-3}$ by Lemma 2.1.
On the other hand, the closed geodesics
$c_1^{2m_1}, c_2^{2m_2}, \ldots c_{\frac{q}{2}}^{2m_{\frac{q}{2}}}$
 contribute $\frac{q}{2}=\frac{n(n+1)d}{4}$ to the summation
$M_{2N+1}+M_{2N+3}+\cdots+M_{2N+nd-3}$ by (\ref{4.22}), (\ref{4.52})
and Lemma 2.1.

Hence
\bea &&\frac{dn^2}{2}=\frac{n(n-1)d}{4}+\frac{n(n+1)d}{4}
\nn\\=&&M_{2N+1}+M_{2N+3}+\cdots+M_{2N+nd-3}
\nn\\=&&b_{2N+1}(\Lambda M)+b_{2N+3}(\Lambda M)+\cdots+b_{2N+nd-3}(\Lambda M)=\frac{dn^2}{2}-1.
\eea
In fact, the second equality holds since $c_k^{2m_k-m}$'s
for $1\le k\le q$ and $m\ge 1$, $c_k^{2m_k}$'s
for $\frac{q}{2}+1\le k\le q$ have no contribution to
$M_{2N+1}+M_{2N+3}+\cdots+M_{2N+nd-3}$ by (\ref{4.51}),
(\ref{4.24}) respectively and Lemma 2.1; $c_k^{2m_k+m}$'s
for $(k,\, m)\notin\Xi$ have no contribution to
$M_{2N+1}+M_{2N+3}+\cdots+M_{2N+nd-3}$
follows from Lemma 2.1 and $i(c_k^{2m_k+m})=2N+i(c_k^m)\ge 2N+nd-1$
by (\ref{4.53}) and (\ref{4.60}) for $m\le \bar m$,
$i(c_k^{2m_k+m})\ge 2N+1+\dim M=2N+nd+1$ by (\ref{4.54}) for $m>\bar m$.
The third equality holds by Step 3 and the last equality
follows from (\ref{4.57}).
This contradiction proves Step 4.

{\bf Step 5.} {\it There exists a closed geodesic $c_{k_\star}$
satisfying $-C(P_{c_{k_\star}})+2\Delta_{k_\star}=nd-1$ among the
$\frac{dn(n+1)}{4}$ closed geodesics found in Step 2.
In particular, $c_{k_\ast}$ and $c_{k_\star}$ are distinct.}

Note that the closed geodesic
$c_{k_\ast}$ found in Step 4 satisfies $k_\ast\in\Omega^{o}_+$.
By the same proof as in Step 4, we find another closed geodesic
$c_{k_\star}$ satisfying $-C(P_{c_{k_\star}})+2\Delta_{k_\star}=nd-1$
with $k_\star\in\Omega^{'o}_+=\Omega^o_-$.
Since $\Omega^o_+\cap\Omega^o_-=\emptyset$,
$c_{k_\ast}$ and $c_{k_\star}$ are distinct.

{\bf Step 6.} {\it Both of  $c_{k_\ast}$ and $c_{k_\star}$ are
irrationally elliptic.}

Firstly we prove

{\bf Claim 2.} {\it The
basic normal form decomposition of the linearized Poincar\'{e} map
$P_{c_{k_\ast}}$  possesses the following
form
\bea f_{c_{k_\ast}}(1)
= R(\th_1)\,\dm\,\cdots\,\dm\,R(\th_{dn-1})\lb{4.64}\eea
with $\frac{\th_{j}}{2\pi}\not\in\Q$ for $1\le j\le dn-1$.
}

By (\ref{4.4}) and (\ref{4.5}), suppose the
basic normal form decomposition of the linearized Poincar\'{e} map
$P_{c_{k_\ast}}$  possesses the following
form
\bea f_{c_{k_\ast}}(1)
&=& R(\th_1)\,\dm\,\cdots\,\dm\,R(\th_r)\,\dm\,D(\lm_1)\dm \,\cdots\,\dm D(\lm_s) \nn\\
& & \dm\,N_2(e^{\sqrt{-1}\aa_{1}},u_{1})\,\dm\,\cdots\,\dm\,N_2(e^{\sqrt{-1}\aa_{r_{\ast}}},u_{r_{\ast}})\,
\nn\\ & &   \dm\,N_2(e^{\sqrt{-1}\bb_{1}},v_{1})\,\dm\,\cdots\,\dm\,N_2(e^{\sqrt{-1}\bb_{r_{0}}},v_{r_{0}}),
\lb{4.65}\eea
where $ N_2(e^{\sqrt{-1}\aa_{j}}, u_j) $'s
are non-trivial and   $N_2(e^{\sqrt{-1}\bb_{j}}, v_j)$'s  are trivial basic normal
forms; $\lm_j\in \{\pm 2\}$ for $1\le j\le s$, $\frac{\th_{j}}{2\pi}\not\in\Q$ for $1\le j\le r$,
$\frac{\aa_{j}}{2\pi}\not\in\Q$ for $1\le j\le r_{\ast}$, $\frac{\bb_{j}}{2\pi}\not\in\Q$ for $1\le j\le r_0$,
and
\be r + s + 2r_{\ast} + 2r_0 = dn - 1. \lb{4.66}\ee
Note that by Theorem 3.5,
\bea C(P_{c_{k_\ast}})&=&\sum_{0<\th<2\pi}S_{P_{c_{k_\ast}}}^-(e^{\sqrt{-1}\th}),
\lb{4.67}\\
\Delta_{k_\ast} &=& \sum_{0<\{m_{k_\ast}\th/\pi\}<\delta}S^-_{P_{c_{k_\ast}}}(e^{\sqrt{-1}\th})\le C(P_{c_{k_\ast}}).
\lb{4.68}\eea
By Lemma 3.1, we have
\bea C(P_{c_{k_\ast}})=\sum_{0<\th<2\pi}S_{P_{c_{k_\ast}}}^-(e^{\sqrt{-1}\th})=r+2r_\ast\le dn-1.
\lb{4.69}\eea
In fact, we have $S^-_{R(\theta_i)}(e^{\sqrt{-1}\th_i})=1$
by (v) of Lemma 3.1,
$S^-_{R(\theta_i)}(e^{\sqrt{-1}\th})=0$ for $\theta\neq\theta_i, 2\pi-\theta_i$
by (\ref{3.5}), $S^-_{R(\theta_i)}(e^{\sqrt{-1}(2\pi-\th_i)})=
S^+_{R(\theta_i)}(e^{\sqrt{-1}\th_i})=0$ by (\ref{3.6}) and (v) of
Lemma 3.1. $S^-_{D(\lambda_i)}(e^{\sqrt{-1}\th})=0$ for
$\theta\in(0, 2\pi)$ by (\ref{3.5}).
$S^-_{N_2(e^{\sqrt{-1}\aa_{j}}, u_j)}(e^{\sqrt{-1}\alpha_j})=1$
by (vi) of Lemma 3.1,
$S^-_{N_2(e^{\sqrt{-1}\aa_{j}}, u_j)}(e^{\sqrt{-1}\th})=0$ for $\theta\neq\alpha_j, 2\pi-\alpha_j$
by (\ref{3.5}), $S^-_{N_2(e^{\sqrt{-1}\aa_{j}}, u_j)}(e^{\sqrt{-1}(2\pi-\alpha_j)})=
S^+_{N_2(e^{\sqrt{-1}\aa_{j}}, u_j)}(e^{\sqrt{-1}\alpha_j})=1$ by (\ref{3.6})
and (vi) of Lemma 3.1.
$S^-_{N_2(e^{\sqrt{-1}\bb_{j}}, v_j)}(e^{\sqrt{-1}\beta_j})=0$
by (vii) of Lemma 3.1,
$S^-_{N_2(e^{\sqrt{-1}\bb_{j}}, v_j)}(e^{\sqrt{-1}\th})=0$ for $\theta\neq\beta_j, 2\pi-\beta_j$
by (\ref{3.5}), $S^-_{N_2(e^{\sqrt{-1}\bb_{j}}, v_j)}(e^{\sqrt{-1}(2\pi-\beta_j)})=
S^+_{N_2(e^{\sqrt{-1}\bb_{j}}, v_j)}(e^{\sqrt{-1}\beta_j})=0$ by (\ref{3.6})
and (vii) of Lemma 3.1. Then (\ref{4.69}) follows from (\ref{3.7}).

By Step 4,
\bea dn-1&=&-C(P_{c_{k_\ast}})+2\Delta_{k_\ast}
\nn\\&\le&-C(P_{c_{k_\ast}})+2C(P_{c_{k_\ast}})
\nn\\&=&C(P_{c_{k_\ast}})=r+2r_\ast\le dn-1.\lb{4.70}\eea
Hence the two inequalities in (\ref{4.70}) must be  equalities,
i.e., we must have
\bea\sum_{0<\{m_{k_\ast}\th/\pi\}<\delta}S^-_{P_{c_{k_\ast}}}(e^{\sqrt{-1}\th})
=\Delta_{k_\ast}=C(P_{c_{k_\ast}})=r+2r_\ast=dn-1.\lb{4.71}\eea
By (\ref{4.66}) and (\ref{4.71}), we have $s=r_0=0$.
Now in order to prove Claim 2, it is sufficient to prove $r_\ast=0$.
Assume the contrary, i.e., there exists at least one
non-trivial $N_2(e^{\sqrt{-1}}\aa_{1},u_{1})$ in (\ref{4.65})
with $\frac{\alpha_1}{2\pi}\notin\Q$.
By (vi) of Lemma 3.1 and (\ref{3.6}),
\bea &&S^-_{N_2(e^{\sqrt{-1}\aa_{1}},u_{1})}(e^{\sqrt{-1}\aa_{1}})=1,
\lb{4.72}\\&&S^-_{N_2(e^{\sqrt{-1}\alpha_1},u_{1})}(e^{\sqrt{-1}(2\pi-\aa_{1})})
=S^+_{N_2(e^{\sqrt{-1}\aa_{1}},u_{1})}(e^{\sqrt{-1}\aa_{1}})=1.
\lb{4.73}\eea
Note that since $\frac{\alpha_1}{\pi}\notin\Q$,
\bea \left\{\frac{m_{k_\ast}(2\pi-\alpha_1)}{\pi}\right\}
=\left\{-\frac{m_{k_\ast}\alpha_1}{\pi}\right\}
=1-\left\{\frac{m_{k_\ast}\alpha_1}{\pi}\right\}.\lb{4.74}
\eea
Thus $\{m_{k_\ast}\alpha_1/\pi\}\in (0, \delta)$ implies
$\{m_{k_\ast}(2\pi-\alpha_1)/\pi\}\notin (0, \delta)$
by choosing $\delta\in (0, \frac{1}{2})$, and then
$$\Delta_{k_\ast}=\sum_{0<\{m_{k_\ast}\th/\pi\}<\delta}S^-_{P_{c_{k_\ast}}}(e^{\sqrt{-1}\th})
<\sum_{0<\th<2\pi}S_{P_{c_{k_\ast}}}^-(e^{\sqrt{-1}\th})=C(P_{c_{k_\ast}})=nd-1.$$
This contradicts to (\ref{4.71}),  hence Claim 2 is true.

Note that by definition, $c_{k_\ast}$ is irrationally elliptic
is equivalent to the following

{\bf Claim 3.} {\it There exist $Q\in\Sp(2dn-2)$
such that  $P_{c_{k_\ast}}=Q^{-1}R(\theta_1)\diamond R(\theta_2)\diamond\cdots\diamond R(\theta_{dn-1})Q$
  with $\frac{\theta_i}{\pi}\notin\Q$
 for $i=1, 2, \ldots,dn-1$ given in (\ref{4.64}).   }

In order to prove claim 3,  we prove the following firstly.

{\bf Claim 4.} {\it If $R(\theta)$ appears in the normal form decomposition
of $P_{c_{k_\ast}}$ in (\ref{4.64}), then $R(2\pi-\theta)$ doesn't  appear.}

Suppose the contrary. Note that  by the same argument as in (\ref{4.74}),  $\{m_{k_\ast}\th/\pi\}\in (0, \delta)$ implies
$\{m_{k_\ast}(2\pi-\th)/\pi\}\notin (0, \delta)$, and then
$\Delta_{k_\ast}<C(P_{c_{k_\ast}})=nd-1$ as above, contradicts to (\ref{4.71}).
This contradiction proves Claim 4.

Now suppose $\omega^{\pm 1}_1,\ldots,\omega^{\pm 1}_{dn-1}$ with
$\omega_i=e^{\sqrt{-1}\theta_i}$
are those eigenvalues of $P_{c_{k_\ast}}$ that satisfy $\frac{\theta_i}{\pi}\notin\Q$,
where $\theta_i$ is given in (\ref{4.64}) for $1\le i\le dn-1$.
Then we have

{\bf Claim 5.}  {\it There exists $Q_1\in \Sp(2dn-2)$ such that
\be Q_1^{-1}P_{c_{k_\ast}}Q_1=R(\vartheta_1)\diamond\cdots\diamond R(\vartheta_{dn-1})
\lb{4.75}\ee
with $\vartheta_i=\theta_i$ or $2\pi-\theta_i$ for $1\le i\le dn-1$.}

In fact, by Theorem 1.6.11 of \cite{Lon3}, we have $X_1\in \Sp(2dn-2)$ such that
\be X_1P_{c_{k_\ast}}X_1^{-1}=S_1\diamond\cdots
\diamond S_{m_1}\diamond S_0,\lb{4.76}\ee
with $S_0\in\Sp(2k_0)$ with $k_0\ge 0$ and $\omega_1\notin\sigma(S_0)$,
$k_i\ge 1$ and $S_i\in\Sp(2k_i)$ is of  the normal form $N_{k_i}(\lambda_i, b_i)$
with $\lambda_i=\omega_1$ or $\omega_1^{-1}$ defined in Section 1.6 of
\cite{Lon3}.

Then by Cases 3 and 4 in Section 1.8 of \cite{Lon3},
if $k_i\ge 3$ for some $1\le i\le m_1$, then $S_i$ can be
connected within $\Omega^0(S_i)$ to $\widetilde{S_i}$ with
$e(\widetilde{S_i})<e(S_i)$. Therefore $e(P_{c_{k_\ast}})<2dn-2$. This contradicts to
 $e(P_{c_{k_\ast}})=2dn-2$ by (\ref{4.64}).
Hence  $S_i=N_{l_i}(\lambda_i, b_i)$ for $l_i=1$ or $2$ and $1\le i\le m_1$.
We next prove that $l_i=1$ for $1\le i\le m_1$. Suppose $l_i=2$ for some $i$.
By (\ref{4.64}), $r_\ast=r_0=0$, thus $S_i$ is not a basic normal form,
and then $S_i\in \mathcal{M}_{\omega_1}^2(4)$. By Case 4 in Section 1.8 of \cite{Lon3},
$S_i$ can be connected within $\Omega^0(S_i)$ to
$R(\vartheta_1)\diamond R(2\pi-\vartheta_1)$. This contradicts to Claim 4.
Now (\ref{4.76}) becomes
\be X_1P_{c_{k_\ast}}X_1^{-1}=R(\vartheta_1)^{\diamond m_1}\diamond S_0,\lb{4.48}\ee
with $\vartheta_1=\theta_1$ or $2\pi-\theta_1$.
We continue the above argument for at most $dn-1$ times and then
obtain Claim 5. Then Claim 3 holds by the uniqueness of the
basic normal form decomposition, cf. Theorem 3.3.
This proves $c_{k_\ast}$ is irrationally elliptic,
the same argument proves $c_{k_\star}$ is irrationally elliptic too.
\hfill\hb

\medskip

\bibliographystyle{abbrv}

\end{document}